# RECOVERING CONVEX BOUNDARIES FROM BLURRED AND NOISY OBSERVATIONS[1]

By Alexander Goldenshluger and Assaf Zeevi

*Haifa University and Columbia University*

We consider the problem of estimating convex boundaries from blurred and noisy observations. In our model, the convolution of an intensity function $f$ is observed with additive Gaussian white noise. The function $f$ is assumed to have convex support $G$ whose boundary is to be recovered. Rather than directly estimating the intensity function, we develop a procedure which is based on estimating the support function of the set $G$. This approach is closely related to the method of *geometric hyperplane probing*, a well-known technique in computer vision applications. We establish bounds that reveal how the estimation accuracy depends on the ill-posedness of the convolution operator and the behavior of the intensity function near the boundary.

**1. Introduction.** We consider the problem of recovering convex boundaries of an image from blurred and noisy observations. Following Hall and Koch [18], we concentrate on the so-called continuous image model (which is explained in what follows). The original image is given by an *intensity function* $f$, supported on a closed convex set $G$ in the $d$-dimensional Euclidean space $\mathbb{R}^d$. In many instances, only a degraded version of the original image is available, where the typical reasons for degradation are *blurring* and *noise*; we refer to [4], pages 51–53, for a detailed discussion of different sources of degradation in imaging. Blurring is modeled using a convolution operation, $K * f$, in which the original image $f$ is subjected to the effects of a *point spread function* $K$. The effects of noise degradation can then be adequately captured by adding a stochastic component to the blurred image.

Received May 2003; revised December 2004.
[1]Supported by the Israel Science Foundation and the Caesarea Rothschild Institute of Haifa University.
*AMS 2000 subject classifications.* 62G05, 62H35.
*Key words and phrases.* Image analysis, convex sets, boundary estimation, deconvolution, support function, geometric probing, rates of convergence.

In this setup, the standard problem is to recover the original image $f$ from its degraded version.

The outlined inverse problem of image deconvolution is the subject of a considerable literature in image processing and statistics. In this paper, however, we focus on directly recovering the *boundary* of the image, rather than estimating the image itself $f$. We note that in virtually all image processing applications, edge detection is one of the standard preliminary analysis steps; this fact serves to motivate our paper. Although many practically important edge detection algorithms for degraded images have been proposed in the image processing literature, their theoretical properties are rarely analyzed. Our goal in this paper is to study the theoretical properties of a particular algorithm which uses geometric probing to estimate the boundary of the image.

Our paper is closely related to two strands of research. The first focuses on nonparametric signal and image deconvolution, which has been extensively studied in statistics; see, for example, [10, 11, 14, 17, 18, 22] and references therein. In this literature, the function $f$ is typically assumed to be smooth, even though its smoothness may be unknown. The second stream of research considers the problem of recovering the boundary of a multidimensional image from *direct* observations; see, for example, [8, 19, 20, 21, 22, 25, 26, 31] and references therein. Recovering boundaries in models that involve indirect observations has been discussed recently in [7] and [15].

The main contributions of this paper are the following. We propose a method for estimating the boundary of the support $G$ of an image $f$, where information on the image is only available in the form of blurred and noisy observations. Our approach focuses on direct recovery of the *support function* of the set $G$; the proposed estimation scheme is akin to the method of *geometric hyperplane probing* which is very common in computer vision applications (see, e.g., [30] and [24]). A similar estimation method has been proposed in [15] for reconstructing convex shapes from noisy observations of their moments. In terms of accuracy of the proposed estimation procedure, we first study the behavior of an estimator of the linear probe functional (see Theorem 1), and subsequently establish upper bounds on accuracy of the support function estimation algorithm (see Theorem 2). We next extend this pointwise result by deriving upper bounds on the *global* estimation accuracy, as measured by the Hausdorff distance between the boundary and its support-function-based estimate (see Theorem 3).

It is worth noting that an alternative approach to the problem studied in this paper is to first deconvolve the function $f$, and subsequently infer the support boundary. However, in our setup the assumptions imposed on the intensity function $f$ are quite weak. In particular, we assume that $f$ is a square integrable bounded function with convex support, with some restrictions imposed on its behavior near the boundary. Under such assumptions, the accuracy of any deconvolution estimator can be arbitrarily bad. Hence, in the absence of additional information on the regularity of $f$, "direct" estimation of the boundary can be advantageous.

The paper is organized in the following way. In Section 2 we formulate the estimation problem, and introduce some preliminaries and notation. Section 3 details our estimation procedure, whose accuracy is analyzed in Section 4. Section 5 contains some discussion, and proofs are all relegated to Section 6.

## 2. Problem formulation.

*The model and observation process.* The image is modeled as an unknown positive function $f \in L_2(\mathbb{R}^d)$, the so-called intensity function, supported on the convex compact set $G$ with nonempty interior in $\mathbb{R}^d$. Suppose that we can observe the process $Y(x)$, given by

$$(1) \qquad dY(x) = (Kf)(x)\,dx + \varepsilon\,dW(x), \qquad x \in \mathbb{R}^d,$$

where $K : \mathcal{D}(K) \mapsto \mathcal{R}(K)$ is a linear transformation with domain $\mathcal{D}(K) \subseteq L_2(\mathbb{R}^d)$, and range $\mathcal{R}(K) \subseteq L_2(\mathbb{R}^d)$, $0 < \varepsilon < 1$, and $(W(x) : x \in \mathbb{R}^d)$ is a standard $d$-dimensional Brownian motion. The observation scheme (1) implies that, for any function $q \in L_2(\mathbb{R}^d)$, we can observe $\langle q, Kf \rangle$ with added zero mean Gaussian noise having variance $\varepsilon^2 \|q\|^2$. Here, $\langle \cdot, \cdot \rangle$ and $\|\cdot\|$ denote the inner product and the corresponding norm in $L_2(\mathbb{R}^d)$ with respect to Lebesgue measure. In what follows, we restrict attention to the case where $K$ is the convolution operator

$$(2) \qquad (Kf)(x) = \int_{\mathbb{R}^d} K(y) f(x-y)\,dy,$$
$$x \in \mathbb{R}^d, f \in L_2(\mathbb{R}^d), K \in L_1(\mathbb{R}^d).$$

The function $K : \mathbb{R}^d \mapsto \mathbb{R}^d$ which models the blurring process is referred to as the *point spread function* in image analysis (cf., e.g., [4], pages 51–53). Following Hall and Koch [18], we refer to (1)–(2) as the *continuous image model*. Our goal is to estimate the boundary of $G$, using observations from (1)–(2).

*Problem geometry and notation.* We introduce some notation that will be used in the sequel. The Euclidean norm of the $d$-vector $x = (x_1, \ldots, x_d)'$ is denoted by $|x|$. Let $B^d = \{x : |x| \leq 1\}$ denote the centered unit ball in $\mathbb{R}^d$; its surface is the unit $d$-sphere $S^{d-1} = \{x : |x| = 1\}$. For a unit vector $u \in S^{d-1}$, we denote by $u^\perp$ the $(d-1)$-dimensional subspace orthogonal to $u$. The centered cube with faces parallel to the coordinate hyperplanes will be designated by $E = \{x : |x_i| \leq 1, 1 \leq i \leq d\}$. In what follows we will consider some simple affine transformations in $\mathbb{R}^d$. In particular, for $u \in S^{d-1}$, let $A_u$ denote the $d \times d$ rotation matrix which maps $(1, 0, \ldots, 0)'$ to the unit vector $u \in S^{d-1}$. For fixed $\tau = (u, r) \in T := S^{d-1} \times [0, 1]$, we define

$$E_\tau := A_u E + (1 + r)u. \tag{3}$$

Thus, $E_\tau$ is obtained by rotating $E$ according to $A_u$, followed by a translation $x \mapsto x + (1 + r)u$. Consequently, $E_\tau$ is the cube centered at $(1 + r)u$ with sides of length 2 and faces parallel to the coordinate hyperplanes of $(u, u^\perp)$.

We assume that $G$, the support of $f$, is a $d$-dimensional convex body, $G \subseteq B^d$. The following concept of a *support function* plays a key role in our estimation procedure.

DEFINITION 1. For $u \in S^{d-1}$, the support function $h(u)$ of the set $G$ is defined by

$$h(u) := \sup\{x'u : x \in G\}.$$

The supporting hyperplane to $G$ with outward normal $u \in S^{d-1}$ is given by $\{x : x'u = h(u)\}$; the support function $h(u)$ at the unit vector $u$ gives the *signed* distance from $o$ to the supporting hyperplane. Throughout the paper we assume that the origin, $o$, belongs to the interior of $G$. This assumption is not restrictive; if it holds, the support function $h(u)$ gives the actual distance from $o$ to the supporting hyperplane.

The key observation that will be exploited in what follows is that there is a one-to-one correspondence between a closed convex set $G$ and its support function, namely,

$$G = \{x : x'u \leq h(u), \forall u \in S^{d-1}\}. \tag{4}$$

For fixed $u \in S^{d-1}$ and small $\eta > 0$, we denote

$$H_u(\eta) := \{x : h(u) - \eta \leq x'u \leq h(u)\}. \tag{5}$$

In words, $H_u(\eta)$ is the set formed by the intersection of $G$ with the half-space $\{x : x'u \geq h(u) - \eta\}$. The volume of the set $H_u(\eta)$ characterizes "massiveness" of $G$ in the direction $u$ in a small vicinity of the support value. It is interesting to note that behavior of the volume of $H_u(\eta)$ as $\eta \to 0$ is related to the tail behavior of the Fourier transform of $\mathbf{1}_G(\cdot)$ in the direction $u \in S^{d-1}$ (see, e.g., the survey paper by Brandolini, Rigoli and Travaglini [5] and references therein). We remark also that $H_u(\eta) = E_{(u, h(u) - \eta)} \cap G$ because $G \subseteq B^d$; this will be used repeatedly in what follows.

**3. The proposed estimation scheme.** Our approach to recovering the boundary of the set $G$ is based on pointwise estimation of its support function $h(\cdot)$. Before developing the details of the proposed approach, let us first highlight the main idea that underlies it. Fix a direction $u$, and imagine that a $d$-dimensional cube is placed outside the unit ball in $\mathbb{R}^d$ so that one of its faces is tangent to the ball in direction $u$. This will be referred to as the *cube probe*. Recall that $G \subseteq B^d$, so by "sliding" this cube toward the origin along the axis $u$, the cube eventually intersects with the boundary of $G$. In the absence of blurring and noise (i.e., $K = I$ and $\varepsilon = 0$), the cube probe will "hit" the boundary when it is exactly at distance $h(u)$ from the origin. As the probe penetrates $G$, it accumulates "mass," where the rate of this mass accumulation is determined by the following: (i) the intensity function ($f$); (ii) the blurring incurred by the convolution operator ($K$); and (iii) the noise level ($\varepsilon$). The main idea is to estimate the distance where the accumulated mass begins to differ significantly from zero, indicating the location of a boundary.

*The probe functional.* For fixed $\tau = (u,r) \in T$, where $T := S^{d-1} \times [0,1]$, let $E_\tau$ be given by (3). We define the *probe functional* by

$$\ell_f(\tau) := \int_{E_\tau} f(x)\,dx = \langle f, \mathbf{1}_{E_\tau} \rangle, \qquad \tau = (u,r) \in T,$$

where $\mathbf{1}_B$ stands for the indicator function of a set $B \in \mathbb{R}^d$. The value of the linear functional $\ell_f(\tau)$ is nothing but the mass accumulated by the cube probe when it is at the location $\tau = (u,r) \in T$. The important property of the probe functional which underlies our construction is the following. For any fixed $u \in S^{d-1}$, if $r \in [h(u), 1]$, then $\ell_f(u,r) = 0$; as $r$ decreases from $h(u)$ to $0$, $\ell_f(u,r)$ grows monotonically, because $f$ is positive. Although we suppose that $f$ is everywhere positive, it is seen from the proofs below that it suffices to assume that $f$ is positive near its support boundary, that is, for $u \in S^{d-1}$, there exists $\Delta > 0$ such that $f(x) \geq 0$ for all $x \in H_u(\Delta)$.

*Estimation of the probe functional.* The first step of our construction is to estimate the probe functional $\ell_f(\tau)$ for given $\tau = (u,r) \in T$. Fix $\delta > 0$, and let $p_\delta$ denote a real-valued infinitely differentiable function $p_\delta \in C^\infty(\mathbb{R}^1)$, equal to 1 on the interval $[-1+\delta, 1-\delta]$ and 0 outside $[-1,1]$. Define

$$(6) \qquad \varphi_\delta(x) := \prod_{j=1}^d p_\delta(x_j), \qquad x = (x_1, \ldots, x_d)'.$$

Clearly, the support of $\varphi_\delta$ is the cube $E = \{x : |x_i| \leq 1, 1 \leq i \leq n\}$. For fixed $\tau = (u,r) \in T$, we define the shifted and rotated version of $\varphi_\delta$ to be

$$(7) \qquad \varphi_{\tau,\delta}(x) := \varphi_\delta(A_u x - (1+r)u).$$

By construction, the support of $\varphi_{\tau,\delta}$ is the shifted and rotated cube, $E_\tau$, defined in (3).

Suppose that $\varphi_{\tau,\delta}$ is in the range of $K^*$, $\mathcal{R}(K^*)$, where $K^*$ is the adjoint operator to $K$. Then, according to the *linear functional strategy* (cf. [1] and [13]), if $\mathcal{D}(K)$ is dense in $L_2(\mathbb{R}^d)$, then there exists $\psi_{\tau,\delta} \in L_2(\mathbb{R}^d)$ such that $\varphi_{\tau,\delta} = K^* \psi_{\tau,\delta}$, and

$$\langle f, \varphi_{\tau,\delta} \rangle = \langle Kf, \psi_{\tau,\delta} \rangle = \langle f, K^* \psi_{\tau,\delta} \rangle \qquad \forall f \in \mathcal{D}(K). \tag{8}$$

Now, for any function $q \in L_1(\mathbb{R}^d) \cap L_2(\mathbb{R}^d)$, let $\widehat{q}$ denote its Fourier transform

$$\widehat{q}(\omega) := \int_{\mathbb{R}^d} q(x) e^{i\omega' x} \, dx, \qquad \omega = (\omega_1, \ldots, \omega_d)'.$$

It is easily verified that, for the convolution operator $K$, the range of the adjoint is given by

$$\mathcal{R}(K^*) = \left\{ q \colon \int_{\mathbb{R}^d} |\widehat{q}(\omega)|^2 |\widehat{K}(\omega)|^{-2} \, d\omega < \infty \right\},$$

and $\psi_{\tau,\delta}$ in (8) is given by

$$\psi_{\tau,\delta}(x) = \frac{1}{(2\pi)^d} \int_{\mathbb{R}^d} \frac{\widehat{\varphi}_{\tau,\delta}(\omega)}{\widehat{K}(-\omega)} e^{-i\omega' x} \, d\omega. \tag{9}$$

The following result states that, under natural assumptions on $K$, $\varphi_{\tau,\delta}$ belongs to $\mathcal{R}(K^*)$ for any $\tau \in T$ and all $\delta > 0$.

ASSUMPTION 1. There exist constants $L > 0$ and $\beta > 0$ such that

$$|\widehat{K}(\omega)| \geq L(1+|\omega|^2)^{-\beta/2} \qquad \forall \omega \in \mathbb{R}^d.$$

LEMMA 1. *Let Assumption 1 hold. Then, $\varphi_{\tau,\delta} \in \mathcal{R}(K^*)$ for any $\tau \in T$ and all $\delta > 0$, and in addition,*

$$\|\psi_{\tau,\delta}\|^2 \leq C_1 L^{-2} \delta^{-2\beta+1}, \tag{10}$$

*where $C_1 = C_1(d, \beta)$ depends on $d$ and $\beta$ only.*

We note that Assumption 1 guarantees the identifiability of $f$ from the observations (1)–(2), and states that the tails of the Fourier transform of the point spread function $K$ cannot decrease to zero faster than the indicated polynomial rate. This assumption is quite standard in deconvolution problems, and corresponds to what is known as a moderately ill-posed problem. We note that the severely ill-posed case, where the tails of $\widehat{K}$ are allowed to decrease at an exponential rate, can be also treated using our methodology. This setup is, however, beyond the scope of our paper.

We are now ready to define the estimate of the probe functional $\ell_f(\tau)$ based on the observations (1). For fixed $\tau = (u,r) \in T$ and $\delta > 0$, we define

$$\tilde{\ell}(\tau;\delta) := \int_{\mathbb{R}^d} \psi_{\tau,\delta}(x)\,dY(x), \tag{11}$$

where $\psi_{\tau,\delta}$ is given in (9). We note that, in view of Lemma 1, the estimator $\tilde{\ell}(\tau;\delta)$ is well defined for all $\tau \in T$ and $\delta > 0$. In what follows, we assume that

$$\sup_{x \in G} |f(x)| \leq M < \infty. \tag{12}$$

The next statement establishes an upper bound on the accuracy of the probe functional estimator $\tilde{\ell}(\tau;\delta)$ for a fixed direction $u \in S^{d-1}$ uniformly over $r \in [0,1]$.

THEOREM 1. *Let Assumption 1 and* (12) *hold, and suppose that $\varepsilon$ is sufficiently small. Define*

$$\delta_* := \left(\frac{\varepsilon}{LM}\sqrt{\ln\frac{1}{\varepsilon}}\right)^{1/(\beta+1/2)}, \tag{13}$$

*and $\tilde{\ell}_*(\tau) := \tilde{\ell}(\tau;\delta_*)$. Then for any fixed $u \in S^{d-1}$, we have*

$$\left\{\mathbb{E}\sup_{r \in [0,1]} |\tilde{\ell}_*(\tau) - \ell_f(\tau)|^2\right\}^{1/2} \leq C_2(d,\beta)\left(M^{\beta-1/2}\frac{\varepsilon}{L}\sqrt{\ln\frac{1}{\varepsilon}}\right)^{1/(\beta+1/2)}. \tag{14}$$

It is important to emphasize that the upper bound of Theorem 1 holds for all square integrable functions $f$ supported on $G \subseteq B^d$ and satisfying (12). If further assumptions on $f$ are introduced, for instance, if $f$ is assumed to be smooth, then the rate of convergence in estimating $\ell_f(\tau)$ may be improved; see further discussion in Section 5.

*The support function estimator.* Based on the estimate $\tilde{\ell}_*(\tau) = \tilde{\ell}_*(u,r)$ of the probe functional $\ell_f(\tau)$, we define the estimator of the support function $h(\cdot)$ at a fixed direction $u \in S^{d-1}$ in the following way:

$$\tilde{h} = \tilde{h}(u) := \max\{r \in [0,1] : \tilde{\ell}_*(u,r) \geq \theta\}, \tag{15}$$

where $\tilde{\ell}_*(u,r) = \tilde{\ell}_*(\tau)$ is given in Theorem 1, and $\theta > 0$ is a threshold to be chosen in the sequel. If $\tilde{\ell}_*(u,r) < \theta$ for all $r \in [0,1]$, we set $\tilde{h}(u) = 0$.

We note that, for $\theta$ sufficiently small, the estimator $\tilde{h}(u)$ is well defined. Indeed, for any fixed $u \in S^{d-1}$, $\tilde{\ell}_*(\tau) = \tilde{\ell}_*(u,r)$, considered as a function of $r$, is a separable Gaussian process on $[0,1]$ that has continuous sample paths with probability 1. The last fact follows from the bounds in the proof of Lemma 2 in Section 6, continuity properties of $\psi_{\tau,\delta}$ as a function of $r$ and the well-known criteria for continuity of sample paths of Gaussian processes (see, e.g., [12], pages 193–194). In particular, this implies that $\tilde{\ell}_*(u, \tilde{h}) = \theta$ with probability 1.

## 4. Accuracy of the support function estimator.

*The class of intensity functions.* To analyze the accuracy of the estimator $\tilde{h}(u)$, we introduce the following class of functions $f$.

DEFINITION 2. We say that a positive function $f \in L_2(\mathbb{R}^d)$ with convex support $G \subseteq B^d$, $o \in \text{int}(G)$, belongs to the class $\mathcal{F}_u(\alpha)$, with $\alpha \geq 1$, if (12) holds and for a fixed direction $u \in S^{d-1}$, there exist $Q > 0$, $\Delta > 0$ such that

$$(16) \qquad \int_{H_u(\eta)} f(x)\,dx \geq Q\eta^\alpha \qquad \forall \eta \in (0, \Delta),$$

where $H_u(\eta)$ is defined in (5).

Several comments on the above definition are in order. The integral on the left-hand side of (16) quantifies the "massiveness" of the image $f$ in the direction $u$. In view of our definition of the probe functional, this integral is nothing but $\ell_f(u, h(u) - \eta)$. Thus, (16) can be equivalently rewritten as

$$(17) \qquad \ell_f(\tau) = \ell_f(u, r) \geq Q|h(u) - r|^\alpha \qquad \forall r \in (h(u) - \Delta, h(u)).$$

It is important to realize that the class $\mathcal{F}_u(\alpha)$ is defined for a fixed direction $u \in S^{d-1}$, and therefore, the constants $Q$, $\Delta$ and $\alpha$ may depend on $u$. Relation (17) imposes restrictions both on the geometry of $G$ in the vicinity of the support value $h(u)$ in direction $u \in S^{d-1}$, and on the behavior of $f$ near the boundary of $G$. The probe functional $\ell_f(u, r)$ quantifies the rate at which the cube probe $E_\tau = E_{(u,r)}$ accumulates mass as it penetrates the set $G$ in a fixed direction $u$ (as $r$ decreases from 1 to 0). The rate of increase in accumulated mass depends on the behavior of the intensity function $f$ near the boundary, and the local curvature of the boundary itself. For example, if $f$ has a discontinuity jump on the boundary, that is, $f = \mathbf{1}_G \tilde{f}$ for $\tilde{f} \geq c > 0$, the probe functional $\ell_f(u, r)$ behaves roughly like the volume of the set $H_u(h(u) - r) = E_\tau \cap G$. For simple geometrical objects, the order of this volume can be easily derived. In general, it turns out that the order of this volume is essentially determined by the (Gauss) curvature of the boundary (see, e.g., [5, 6]). We note that the parameter $\alpha$ indexing the class $\mathcal{F}_u(\alpha)$ satisfies $\alpha \geq 1$ whenever $G$ is a convex set. The next examples illustrate how the curvature of the boundary of $G$ and the behavior of $f$ near the boundary determine the index $\alpha$. In what follows we use the term *sharp boundary* when the intensity function $f$ is discontinuous on the boundary of the support set.

EXAMPLE 1. Assume that $d=2$ and consider the case of a sharp boundary. Let $G$ be a convex polygon, and for simplicity, suppose that $f = \mathbf{1}_G$. Then for any direction which is not perpendicular to the side of the polygon, we have $\alpha = 2$. In this case the corresponding supporting line intersects the boundary of $G$ in a vertex of the polygon, and the constant $Q$ depends in an obvious way on the angle between two adjacent sides of the vertex. If the direction $u \in S^1$ is perpendicular to the side of the polygon, then the supporting line contains that side; in this case $\alpha = 1$.

EXAMPLE 2. We again consider the case of $d=2$, and first assume that the boundary is *sharp*. If $G$ is a circle or an ellipse, then $\alpha = 3/2$ for any direction. It turns out that this case is rather general. Let $x_u$ denote the point on the boundary $\partial G$ where the support value in direction $u$ is attained. If $\partial G$ has nonzero curvature at $x_u$, then $\alpha = 3/2$ for that direction. Now assume that the boundary is *nonsharp*. Specifically, without loss of generality, we let $u = (0,1) \in S^1$, and suppose that $\partial G$ can be represented as $x_2 = g(x_1) := -x_1^2 + c$, with $c > 0$, in the vicinity of the origin. This corresponds to the case of nonzero boundary curvature in direction $u$. Let $f(x_1, x_2) \geq |x_2 - g(x_1)|^\gamma$ for some $\gamma \geq 0$ ($\gamma = 0$ corresponds to the sharp boundary case). Then it is easily verified that $\alpha = \gamma + 3/2$. If the curvature of $\partial G$ vanishes at $x_u$, the exponent $\alpha$ differs from $\gamma + 3/2$. For instance, if $\partial G$ is represented as the graph of the function $x_2 = -x_1^{2q} + c$, $q \geq 1$, then $\alpha = \gamma + 1 + (2q)^{-1}$. We note, however, that if the boundary is smooth, then the set of points on the boundary where the curvature vanishes has zero Lebesgue measure.

EXAMPLE 3. In $\mathbb{R}^d$ we consider the case of a sharp boundary. If $G$ is an ellipsoid in $\mathbb{R}^d$, then $\alpha = (d+1)/2$ for any direction $u \in S^{d-1}$. Similarly to the previous example, if for some $u \in S^{d-1}$ the Gauss curvature of $\partial G$ does not vanish at $x_u$, then $\alpha = (d+1)/2$ as well. If $G$ is a parallelotope, $\alpha$ takes values in the set $\{1, 2, \ldots, d\}$ depending on the direction $u$. In particular, if $u$ is perpendicular to a face of the parallelotope, then $\alpha = 1$. Without loss of generality, consider now the situation when $u = (0, \ldots, 0, 1) \in S^{d-1}$ and $\partial G$ can be represented as $x_d = g(x_1, \ldots, x_{d-1}) := -x_1^{2q_1} \cdots x_{d-1}^{2q_{d-1}} + c$ for some $q_i \geq 1$, $i = 1, \ldots, d-1$, and $c > 0$. Then it is not difficult to verify that $\alpha = 1 + \sum_{j=1}^{d-1} (2q_j)^{-1}$.

*Accuracy bounds.* We are now ready to establish upper bounds on the accuracy of the proposed estimation procedure.

THEOREM 2. *Let Assumption* 1 *hold. Fix an arbitrary* $u \in S^{d-1}$, *and let* $\tilde{h}_*(u)$ *be the estimator of the support function given by* (15), *with* $\theta = \theta_*$

*taken to be*

$$\theta_* := \left(\frac{\varepsilon}{L} M^{\beta-1/2} \sqrt{C_3 \ln \frac{1}{\varepsilon}}\right)^{1/(\beta+1/2)}, \qquad C_3 > 0. \tag{18}$$

*Then, for sufficiently small $\varepsilon$,*

$$\sup_{f \in \mathcal{F}_u(\alpha)} \{\mathbb{E}|\tilde{h}_*(u) - h(u)|^2\}^{1/2} \leq C_4 Q^{-1/\alpha} \left\{M^{\beta-1/2} \frac{\varepsilon}{L} \sqrt{\ln \frac{1}{\varepsilon}}\right\}^{1/(\alpha(\beta+1/2))},$$

*where $C_3$ and $C_4$ may depend on $d$, $\beta$ and $\Delta$ only.*

We note that since the blurring kernel $K$ is assumed to be known, the dependence of the above threshold $\theta_*$ on $L$ and $\beta$, which characterize the behavior of $K$, is not restrictive.

Theorem 2 gives an upper bound on the *pointwise risk* of the proposed estimator. Based on the estimate $\tilde{h}_*$, we can define the global estimator of the set $G$ as

$$\tilde{G}_* := \{x : x'u \leq \tilde{h}_*(u), \forall u \in S^{d-1}\}. \tag{19}$$

We note that $\tilde{G}_*$ is a convex set, by construction, because it is given by the intersection of half-spaces formed by the estimated supporting hyperplanes. In order to measure global accuracy of $\tilde{G}_*$ as an estimator of $G$, it is natural to use the following family of $L_p$-metrics on the class of all convex bodies in $\mathbb{R}^d$. If $G_1$ and $G_2$ are convex bodies with support functions $h_{G_1}$ and $h_{G_2}$, respectively, then the $L_p$-metric is defined as

$$\Delta_p(G_1, G_2) := \left\{\int_{S^{d-1}} |h_{G_1}(u) - h_{G_2}(u)|^p \, d\lambda_{d-1}(u)\right\}^{1/p}, \qquad 1 \leq p \leq \infty,$$

where $\lambda_{d-1}$ is the spherical Lebesgue measure on $S^{d-1}$. It is remarkable that the $\Delta_\infty$-metric is nothing but the usual Hausdorff distance between the sets (see, e.g., [29], Section 1.8). We note also that the metric $\Delta_p$ may be viewed as the general $L_p$-metric of [3] specialized to convex sets; the cited paper argues that such a metric has attractive properties for image processing applications.

To state the upper bound on the risk of $\tilde{G}_*$ we first introduce the following "global" version of the functional class $\mathcal{F}_u(\alpha)$.

DEFINITION 3. We say that a positive function $f$ with convex support $G \subseteq B^d$, $o \in \text{int}(G)$, belongs to the class $\mathcal{F}(\alpha)$, if (12) holds and for *all* $u \in S^{d-1}$, there exist $Q > 0$, $\Delta > 0$ and $\alpha \geq 1$ such that

$$\int_{H_u(\eta)} f(x) \, dx \geq Q\eta^\alpha \qquad \forall \eta \in (0, \Delta),$$

where $H_u(\eta)$ is defined in (5).

In contrast to Definition 2, here we require that (16) holds for all directions $u \in S^{d-1}$ with the same constants $Q$, $\alpha$ and $\Delta$. This class may be adequate, for example, for purposes of describing indicator functions of convex sets with smooth boundaries having everywhere nonvanishing Gauss curvature.

It follows immediately from Theorem 2 that

$$\sup_{f \in \mathcal{F}(\alpha)} \{\mathbb{E}\,\Delta_2^2(\tilde{G}_*, G)\}^{1/2} \leq C_5 Q^{-1/\alpha} \left\{ M^{\beta-1/2} \frac{\varepsilon}{L} \sqrt{\ln \frac{1}{\varepsilon}} \right\}^{1/(\alpha(\beta+1/2))}.$$

Our next result concerns global estimation accuracy as measured by the Hausdorff distance.

THEOREM 3. *Let $\tilde{G}_*$ be given by* (19). *Then, under the conditions of Theorem* 2,

$$\sup_{f \in \mathcal{F}(\alpha)} \{\mathbb{E}\Delta_\infty^2(\tilde{G}_*, G)\}^{1/2} \leq C_6 Q^{-1/\alpha} \left\{ M^{\beta-1/2} L^{-1} \varepsilon \sqrt{\ln \frac{1}{\varepsilon}} \right\}^{1/(\alpha(\beta+1/2))}.$$

**5. Discussion.** We now turn to a few brief comments on the main results of this paper.

1. The proposed algorithm is based on direct estimation of the edge by means of the geometric hyperplane probing scheme. As our results show, the accuracy of this method is determined by the ill-posedness of the convolution operator $K$ (given by the value of $\beta$) and the "massiveness" of the image in the $u$-direction (near the boundary) as quantified by the index $\alpha$. The latter is determined by the local behavior of the intensity function $f$ near the boundary and the local geometrical properties of the support $G$ in the vicinity of the estimated support value; these two factors determine the index $\alpha$ in the manner illustrated in Examples 1–3 in Section 4.

2. The functional class $\mathcal{F}_u(\alpha)$ involves only weak assumptions on the behavior of the intensity function. Specifically, $f \in \mathcal{F}_u(\alpha)$ is only assumed to be bounded and square integrable with convex support, satisfying certain "massiveness" conditions near the boundary. We would like to emphasize that, under these conditions uniform bounds on the accuracy of a deconvolution-based estimator can be arbitrarily bad. This is due to the fact that the functional class $\mathcal{F}_u(\alpha)$ is too broad. Hence, our approach, which is built on "directly" estimating the boundary, can be applied in instances where the deconvolution-based estimator is not appropriate.

3. The basic ingredient in our construction is the estimation of the probe functional corresponding to the hyperplane probing scheme. Using Theorem 1 and results in [9], we can show that our estimator of the probe functional is rate optimal (up to a logarithmic term). The proposed boundary estimator inherits this convergence rate.

4. If additional regularity conditions on $f$ are imposed, then one can improve the rates given in the main results of this paper. In particular, there is a close connection between the setup and problem formulation in our paper, and that of estimating a change-point from indirect observations in the one-dimensional case (see, e.g., [27]). The recent paper by Goldenshluger, Tsybakov and Zeevi [16] studies the problem of change-point estimation in a function $f$ in the white noise convolution setup. In that paper it is shown that minimax rates of convergence are determined by smoothness of $f$ away from the change-point, and by the ill-posedness of the convolution operator. Specifically, further smoothness of the function $f$ away from the change-point results in better accuracy in estimating the change-point. (This stands in contrast to the case of direct observations.) This suggests that the optimal rate of convergence in estimating convex boundaries from indirect observations depends on the smoothness of the intensity function $f$. Establishing this result rigorously remains an open problem.

## 6. Proofs.

PROOF OF LEMMA 1. We assume that $\delta > 0$ is fixed, and note that, by (7),

$$
\begin{aligned}
\widehat{\varphi}_{\tau,\delta}(\omega) &= \int_{\mathbb{R}^d} \varphi_{\tau,\delta}(x) e^{i\omega'x} \, dx \\
&= \int_{\mathbb{R}^d} \varphi_\delta(A_u x - (1+r)u) e^{i\omega'x} \, dx \\
&= \exp\{i(1+r)u' A_u \omega\} \int_{\mathbb{R}^d} \varphi_\delta(y) e^{i\omega' A_u' y} \, dy \\
&= \exp\{i(1+r)u' A_u \omega\} \widehat{\varphi}_\delta(A_u \omega).
\end{aligned}
\tag{20}
$$

It follows from (6) that $\widehat{\varphi}_\delta(\omega) = \prod_{j=1}^d \widehat{p}_\delta(\omega_j)$ for all $\omega = (\omega_1, \ldots, \omega_d)' \in \mathbb{R}^d$. The function $\widehat{p}_\delta(\cdot)$ is rapidly decreasing, that is, $|\widehat{p}_\delta(\lambda)| \leq C_k(1+|\lambda|)^{-k}$ for all $k$. Using this fact, and taking into account Assumption 1 and (20), we conclude that $\widehat{\varphi}_{\tau,\delta}(\cdot)/\widehat{K}(-\cdot) \in L_1(\mathbb{R}^d) \cap L_2(\mathbb{R}^d)$ so that $\psi_{\tau,\delta}(\cdot)$ is well defined in (9). This also implies that $\varphi_{\tau,\delta} \in \mathcal{R}(K^*)$.

Now we prove (10). Straightforward algebra shows that

$$
\begin{aligned}
\|\psi_{\tau,\delta}\|^2 &= \int_{\mathbb{R}^d} \left|\frac{\widehat{\varphi}_\delta(A_u \omega)}{\widehat{K}(-\omega)}\right|^2 d\omega \\
&\leq L^{-2} \int_{\mathbb{R}^d} |\widehat{\varphi}_\delta(A_u \omega)|^2 (1+|\omega|^2)^\beta \, d\omega \\
&= L^{-2} \int_{\mathbb{R}^d} |\widehat{\varphi}_\delta(\omega)|^2 (1+|\omega|^2)^\beta \, d\omega \\
&= L^{-2} \int_{\mathbb{R}^d} \prod_{j=1}^d |\widehat{p}_\delta(\omega_j)|^2 (1+|\omega|^2)^\beta \, d\omega \\
&\leq c_1 L^{-2} \int_{\mathbb{R}^1} |\widehat{p}_\delta(\lambda)|^2 |\lambda|^{2\beta} \, d\lambda,
\end{aligned}
\tag{21}
$$

where $c_1 = c_1(\beta, d)$ is the constant depending on $d$ and $\beta$ only. Thus, it is sufficient to bound the last integral in (21). To this end, we recall the standard way to construct the function $p_\delta \in C^\infty(\mathbb{R}^1)$ with aforementioned properties (see, e.g., [23] and [28]). Let $\gamma(x) = e^{-1/[x(1-x)]}$. Then, on the interval $[-1, -1+\delta]$, where $p_\delta$ climbs from zero to one, set

$$p_\delta(x) = \frac{c_2}{\delta} \int_{-1}^{x} \gamma\left(\frac{y+1}{\delta}\right) dy, \qquad x \in [-1, -1+\delta],$$

where $c_2$ is an absolute constant. First assume that $\beta$ is an integer; then

$$p_\delta^{(\beta)}(x) = \frac{c_2}{\delta^\beta} \gamma^{(\beta-1)}\left(\frac{x+1}{\delta}\right), \qquad x \in [-1, -1+\delta],$$

and therefore

$$\int_{\mathbb{R}^1} |p_\delta^{(\beta)}(x)|^2 \, dx = \frac{c_2^2}{\delta^{2\beta}} \int_{\mathbb{R}^1} \left|\gamma^{(\beta-1)}\left(\frac{x+1}{\delta}\right)\right|^2 dx$$

$$= \frac{c_2^2}{\delta^{2\beta-1}} \int_{\mathbb{R}^1} |\gamma^{(\beta-1)}(x)|^2 \, dx$$

$$\leq \frac{c_3}{\delta^{2\beta-1}},$$

where $c_3 = c_3(\beta)$ depends on $\beta$ only. Combining this with (21), we obtain (10) for integer $\beta$. For general $\beta$, the result follows from the standard interpolation inequalities for Sobolev spaces (cf. [2], page 127). □

We now provide an auxiliary lemma that will be used repeatedly in the proofs below. Define

(22) $$X(\tau) = X(u, r) := \varepsilon \int_{\mathbb{R}^d} \psi_{\tau, \delta_*}(x) \, dW(x), \qquad r \in [0, 1],$$

where $\tau = (u, r) \in T$, $\psi_{\tau, \delta_*}$ is given by (9) and $\delta_*$ is defined in (13). In the sequel we will treat $X(\tau)$ as a random process over the index set $T = S^{d-1} \times [0, 1]$. When $u \in S^{d-1}$ is fixed, we have the random process $\{X(u, r) : r \in [0, 1]\}$. This distinction will always be clear from the context.

Obviously, $\{X(\tau)\}$ is a zero mean Gaussian process. Let

$$\sigma_X^2 := \sup_{\tau \in T} \mathbb{E}|X(\tau)|^2,$$

and note that, in view of Lemma 1 and (13), we have that

(23) $$\begin{aligned} \sigma_X^2 &= \sup_{\tau \in T} \varepsilon^2 \|\psi_{\tau, \delta_*}\|^2 \\ &\leq \frac{c_1 \varepsilon^2}{L^2 \delta_*^{2\beta-1}} \\ &\leq c_2 \left(M^{2\beta-1} \frac{\varepsilon^2}{L^2}\right)^{1/(\beta+1/2)}. \end{aligned}$$

LEMMA 2. *Let Assumption* 1 *and condition* (12) *hold.*

(i) *For any fixed* $u \in S^{d-1}$ *and for all* $\theta > 2\sigma_X$, *we have that*

$$\mathbb{P}\bigg\{\sup_{r\in[0,1]}|X(u,r)| \geq \theta\bigg\} \leq \frac{C_1 M \theta}{\sigma_X^2}\exp\bigg\{-\frac{\theta^2}{2\sigma_X^2}\bigg\},$$

*where* $C_1 = C_1(d,\beta)$ *depends on* $d$ *and* $\beta$ *only.*

(ii) *For all* $\theta > 2\sigma_X$, *we have that*

$$\mathbb{P}\bigg\{\sup_{\tau\in T}|X(\tau)| \geq \theta\bigg\} \leq C_2\bigg(\frac{\theta}{\sigma_X^2}\bigg)^d \exp\bigg\{-\frac{\theta^2}{2\sigma_X^2}\bigg\},$$

*where* $C_2$ *may depend on* $d$, $\beta$, $M$ *and* $L$ *only.*

PROOF. (i) Define the semi-metric

$$\mu(\rho,r) := (\mathbb{E}|X(u,\rho) - X(u,r)|^2)^{1/2}.$$

To prove the inequality in the lemma, we verify the conditions of Proposition A.2.7 in [32]. Specifically, we need to evaluate the minimal number of balls of radius $\nu$ with respect to the semi-metric $\mu$ covering the index set $[0,1]$. We have that

$$\begin{aligned}
\mu^2(\rho,r) &= \varepsilon^2 \|\psi_{(u,\rho),\delta_*} - \psi_{(u,r),\delta_*}\|_2^2 \\
&= \varepsilon^2 \int_{\mathbb{R}^d} \bigg|\frac{\widehat{\varphi}_{(u,\rho),\delta_*}(\omega) - \widehat{\varphi}_{(u,r),\delta_*}(\omega)}{\widehat{K}(-\omega)}\bigg|^2 dw \\
&\leq \varepsilon^2 \int_{\mathbb{R}^d} |\widehat{\varphi}_{(u,\rho),\delta_*}(\omega) - \widehat{\varphi}_{(u,r),\delta_*}(\omega)|^2 (1+|\omega|^2)^{\beta} d\omega \\
&\stackrel{(a)}{\leq} \varepsilon^2 \int_{\mathbb{R}^d} |\widehat{\varphi}_{\delta_*}(A_u\omega)|^2 (1+|\omega|^2)^{\beta} |1 - e^{i(\rho-r)u'A_u\omega}|^2 d\omega \\
&\stackrel{(b)}{\leq} \varepsilon^2 |\rho-r|^2 \int_{\mathbb{R}^d} |\widehat{\varphi}_{\delta_*}(\omega)|^2 (1+|\omega|^2)^{\beta} |u'\omega|^2 d\omega \\
&\leq \varepsilon^2 |\rho-r|^2 \int_{\mathbb{R}^d} |\widehat{\varphi}_{\delta_*}(\omega)|^2 (1+|\omega|^2)^{\beta+1} d\omega \\
&\stackrel{(c)}{\leq} c_1 |\rho-r|^2 \frac{\varepsilon^2}{L^2 \delta_*^{2\beta+1}} \\
&= c_1 M^2 |\rho-r|^2,
\end{aligned}$$

where (a) follows from (20); (b) follows from a change of variable noting that $A_u$ is a rotation matrix, and the standard bound $|1 - \exp\{ix\}| \leq |x|$; (c) follows from the same argument as in the proof of Lemma 1; and the last equality follows by definition of $\delta_*$ given in (13). Therefore, the minimal number of balls of radius $\nu$, with respect to the semi-metric $\mu$, that covers the interval $[0,1]$ does not exceed $N(\nu; \mu, [0,1]) = c_2 M \nu^{-1}$. Applying Proposition A.2.7 from [32] with $\varepsilon_0 = \sigma_X$, we obtain that, for some constant $c_3 > 0$,

$$\mathbb{P}\left\{\sup_{r \in [0,1]} |X(u,r)| \geq \theta\right\} \leq \frac{c_3 M \theta}{\sigma_X^2} \exp\left\{-\frac{\theta^2}{2\sigma_X^2}\right\},$$

provided that $\theta > 2\sigma_X$. This concludes the proof of (i).

(ii) Now we consider

$$\mu(\tau, t) = (\mathbb{E}|X(\tau) - X(t)|^2)^{1/2},$$

where $\tau = (u, \rho)$, $t = (v, r)$ and $\tau, t \in T = S^{d-1} \times [0,1]$. Then, similarly to the above considerations, we have

$$\mu^2(\tau, t) \leq \varepsilon^2 \int_{\mathbb{R}^d} |\widehat{\varphi}_{\tau, \delta_*}(\omega) - \widehat{\varphi}_{t, \delta_*}(\omega)|^2 (1 + |\omega|^2)^\beta \, d\omega$$

$$\leq c_4 \varepsilon^2 \left\{ \int_{\mathbb{R}^d} |\widehat{\varphi}_{\delta_*}(A_u \omega)|^2 |e^{i(1+r)u' A_u \omega} - e^{i(1+\rho)v' A_v \omega}|^2 (1 + |\omega|^2)^\beta \, d\omega \right.$$

$$\left. + \int_{\mathbb{R}^d} |\widehat{\varphi}_{\delta_*}(A_u \omega) - \widehat{\varphi}_{\delta_*}(A_v \omega)|^2 (1 + |\omega|^2)^\beta \, d\omega \right\}$$

$$=: c_4 \varepsilon^2 \{J_1 + J_2\}.$$

Straightforward algebra shows that

$$\begin{aligned}
J_1 &\leq c_5 \int_{\mathbb{R}^d} |\widehat{\varphi}_{\delta_*}(A_u \omega)|^2 |(1+r)u' A_u \omega - (1+\rho)v' A_v \omega|^2 (1 + |\omega|^2)^\beta \, d\omega \\
&= c_5 (1 + \rho)^2 \\
&\quad \times \int_{\mathbb{R}^d} |\widehat{\varphi}_{\delta_*}(A_u \omega)|^2 \left\{ \omega' A_u'(u - v) + \omega'(A_u - A_v)'v - \frac{r - \rho}{1 + \rho} \omega' A_v' v \right\}^2 \\
&\quad \times (1 + |\omega|^2)^\beta \, d\omega \\
&\stackrel{(a)}{\leq} c_6 \int_{\mathbb{R}^d} |\widehat{\varphi}_{\delta_*}(A_u \omega)|^2 |\omega|^2 \{|u - v|^2 + (r - \rho)^2\} (1 + |\omega|^2)^\beta \, d\omega \\
&\leq c_7 |t - \tau|^2 \int_{\mathbb{R}^d} |\widehat{\varphi}_{\delta_*}(A_u \omega)|^2 (1 + |\omega|^2)^{\beta+1} \, d\omega \\
&\stackrel{(b)}{\leq} c_8 |t - \tau|^2 (L \delta_*^{2\beta+1})^{-1},
\end{aligned}$$

(24)

where (a) follows from the fact that $A_u$ and $A_v$ are orthogonal matrices, hence, from a singular value decomposition, one has that the spectral norm of the difference between $A_u$ and $A_v$ is $O(|u - v|)$ when $u, v \in S^{d-1}$ are close; and (b) is obtained as in (i).

In order to bound $J_2$ from above, we first assume that $\beta$ is an integer. We introduce the standard multi-index notation: $D^k = D_1^{k_1} \cdots D_d^{k_d} =$ $\partial^{|k|}/\partial x_1^{k_1}\cdots\partial x_d^{k_d}$ is the differential operator of order $|k| = k_1 + \cdots + k_d$, $k = (k_1, \ldots, k_d)$, and for $x = (x_1, \ldots, x_d)$, we write $x^k = x_1^{k_1}\cdots x_d^{k_d}$. Further, we note that $\widehat{\varphi}_\delta(\cdot)$ is infinitely differentiable; hence, there exists $0 \leq \zeta \leq 1$ such that

$$(25) \quad \widehat{\varphi}_{\delta_*}(A_u\omega) = \widehat{\varphi}_{\delta_*}(A_v\omega) + (A_v\omega - A_u\omega)'\nabla\widehat{\varphi}_{\delta_*}(A_v\omega + \zeta A_u\omega),$$

where $\nabla := (D_1, \ldots, D_p)$ is the gradient vector. Using (25), we can write

$$(26) \quad J_2 \leq c_8|u-v|^2 \int_{\mathbb{R}^d} |\nabla\widehat{\varphi}_{\delta_*}((A_v + \zeta A_u)\omega)|^2 (1 + |\omega|^2)^{\beta+1}\, d\omega,$$

and our current goal is to bound the last integral. For $u$ and $v$ close to each other, the matrix $I + \zeta A_v' A_u$ is nonsingular, thus, using a change-of-variables in the above integral, we see that it suffices to bound from above the integral $J_3 := \int_{\mathbb{R}^d} |\nabla\widehat{\varphi}_{\delta_*}(\omega)|^2 (1 + |\omega|^2)^{\beta+1}\, d\omega$. To this end, we note that

$$J_3 \leq c_9 \sum_{|k|\leq\beta+1} \int_{\mathbb{R}^d} |D^k\{(-ix)\varphi_{\delta_*}(x)\}|^2\, dx$$

$$\leq c_{10} \int_{-1}^{1} |p_\delta^{(\beta+1)}(x)|^2\, dx \leq c_{11}\delta_*^{-(2\beta+1)},$$

where the last two inequalities follow from the definition of $\varphi_{\delta_*}$ (and $p_{\delta_*}$) and from the same argument as in the proof of Lemma 1. Combining this with (24) and (26), and substituting (13) for $\delta_*$, we obtain $\mu^2(t,\tau) \leq c_{12}|t-\tau|^2$. Thus, the minimal number of balls of the radius $\nu$ in the semi-metric $\mu$ covering the index set $T = S^{d-1} \times [0,1]$ does not exceed $N(\nu;\mu,T) \leq c_{13}\nu^{-d}$. Therefore, applying Proposition A.2.7 from [32], we come to the result (ii). The statement of the lemma for noninteger $\beta$ follows from standard interpolation inequalities for Sobolev spaces (cf. [2], page 127). $\square$

PROOF OF THEOREM 1. It follows from (1), (11) and (8) that

$$\tilde{\ell}(\tau;\delta) = \langle\psi_{\tau,\delta}, Kf\rangle + \varepsilon\int_{\mathbb{R}^d}\psi_{\tau,\delta}(x)\, dW(x)$$

$$= \langle f, \varphi_{\tau,\delta}\rangle + \varepsilon\int_{\mathbb{R}^d}\psi_{\tau,\delta}(x)\, dW(x).$$

By definition $\ell_f(\tau) = \langle f, \mathbf{1}_{E_\tau}\rangle$. Therefore,

$$(27) \quad \begin{aligned} |\tilde{\ell}(\tau;\delta) - \ell_f(\tau)|^2 &\leq 2|\langle f, \varphi_{\tau,\delta} - \mathbf{1}_{E_\tau}\rangle|^2 + 2\left|\varepsilon\int_{\mathbb{R}^d}\psi_{\tau,\delta}(x)\, dW(x)\right|^2 \\ &\leq c_1(d)M^2\delta^2 + 2|X(u,r)|^2, \end{aligned}$$

where we have taken into account the definition of $\varphi_{\tau,\delta}(\cdot)$ and (12); $X(u,r)$ is defined in (22). In view of Lemma 2, we have, for any $\bar{t} > 2\sigma_X$,

$$\begin{aligned}
\mathbb{E} \sup_{r \in [0,1]} |X(u,r)|^2 &= \int_0^\infty 2t \mathbb{P}\left\{ \sup_{r \in [0,1]} |X(u,r)| > t \right\} dt \\
&\leq 2\bar{t}^2 + c_3 \int_{\bar{t}}^\infty t^2 \sigma_X^{-2} \exp\{-t^2/(2\sigma_X^2)\} dt \\
&= 2\bar{t}^2 + c_4 \sigma_X \int_{\bar{t}^2/2\sigma_X^2} \sqrt{s} \exp\{-s\} ds \\
&\leq 2\bar{t}^2 + c_4 \sigma_X \exp\{-\bar{t}^2/(4\sigma_X^2)\},
\end{aligned}$$

where the last step uses Jensen's inequality. Minimizing with respect to $\bar{t}$, we obtain $\mathbb{E} \sup_{r \in [0,1]} |X(u,r)|^2 \leq c_5 \sigma_X^2 \ln(1/\sigma_X)$, and the statement of the theorem follows from (23). $\square$

PROOF OF THEOREM 2. Throughout the proof we fix a direction $u \in S^{d-1}$ and therefore, for brevity, omit the argument $u$. With slight abuse of notation, we write $h$ for $h(u)$, $\ell_f(r)$ for $\ell_f(u,r)$ and so on. We denote by $c_1, c_2, \ldots$ positive constants depending on $d$, $\beta$ and $\Delta$ only.

For the estimate $\tilde{h}_*$, we have

$$\begin{aligned}
\mathbb{E}|\tilde{h}_* - h|^2 =\ & \mathbb{E}[|\tilde{h}_* - h|^2 \mathbf{1}\{h - \Delta \leq \tilde{h}_* \leq h\}] \\
& + \mathbb{E}[|\tilde{h}_* - h|^2 \mathbf{1}\{\tilde{h}_* \leq h - \Delta\}] \\
& + \mathbb{E}[|\tilde{h}_* - h|^2 \mathbf{1}\{\tilde{h}_* > h\}] \\
=:\ & I_1 + I_2 + I_3,
\end{aligned}$$
(28)

where $\Delta$ is defined in (17). We bound $I_1$, $I_2$ and $I_3$ separately.

In view of (17), we have

$$\begin{aligned}
|\tilde{h}_* - h|^2 \mathbf{1}\{h - \Delta \leq \tilde{h}_* \leq h\} &\leq Q^{-2/\alpha} |\ell_f(\tilde{h}_*)|^{2/\alpha} \\
&\leq Q^{-2/\alpha} \{|\ell_f(\tilde{h}_*) - \tilde{\ell}_*(\tilde{h}_*)| + |\tilde{\ell}_*(\tilde{h}_*)|\}^{2/\alpha} \\
&\leq c_1 Q^{-2/\alpha} \{|\ell_f(\tilde{h}_*) - \tilde{\ell}_*(\tilde{h}_*)|^{2/\alpha} + |\tilde{\ell}_*(\tilde{h}_*)|^{2/\alpha}\} \\
&\leq c_1 Q^{2/\alpha} \left\{ \sup_{r \in [h-\Delta, h]} |\tilde{\ell}_*(r) - \ell_f(r)|^{2/\alpha} + \theta_*^{2/\alpha} \right\},
\end{aligned}$$
(29)

where we have used the fact that $\tilde{\ell}_*(\tilde{h}_*) = \theta_*$ with probability 1. Taking the expectation and using Theorem 1, we obtain

$$I_1 \leq c_2 Q^{-2/\alpha} \left\{ \left( M^{\beta-1/2} \frac{\varepsilon}{L} \sqrt{\ln \frac{1}{\varepsilon}} \right)^{2/(\alpha(\beta+1/2))} + \theta_*^{2/\alpha} \right\}.$$

For the second term $I_2$, we have

$$I_2 = \mathbb{E}[|\tilde{h}_* - h|^2 \mathbf{1}\{\tilde{h}_* < h - \Delta\}]$$

$$
\begin{aligned}
&\leq \mathbb{P}\{\tilde{h}_* < h - \Delta\} \\
&\leq \mathbb{P}\{\tilde{\ell}_*(h-\Delta) \leq \theta_*\} \\
&= \mathbb{P}\{\ell_f(h-\Delta) + \tilde{\ell}_*(h-\Delta) - \ell_f(h-\Delta) \leq \theta_*\} \\
&\leq \mathbb{P}\{|\ell_f(h-\Delta) - \tilde{\ell}_*(h-\Delta)| \geq \ell_f(h-\Delta) - \theta_*\} \\
&\leq \mathbb{P}\{|X(u, h-\Delta)| \geq \ell_f(h-\Delta)/2\} \\
&\stackrel{(a)}{\leq} \mathbb{P}\{|\mathcal{N}(0,1)| \geq Q\Delta^\alpha (\varepsilon \|\psi_{\tau,\delta_*}\|)^{-1}/2\} \\
&\stackrel{(b)}{\leq} c_3 \exp\{-c_4 \varepsilon^{-2/(\beta+1/2)}\},
\end{aligned}
$$

where (a) follows from definition of $X(u,r)$, the choice of $\delta_*$ and the fact that $\varepsilon$ is sufficiently small, and (b) follows from Lemma 1. Thus, $I_2 = o(I_1)$ as $\varepsilon \to 0$.

It remains to bound the third error term $I_3$. We have

$$
\begin{aligned}
I_3 = \mathbb{E}[|\tilde{h}_* - h|^2 \mathbf{1}\{\tilde{h}_* > h\}] &\leq \mathbb{P}\{\tilde{h}_* > h\} \\
&= \mathbb{P}\Big\{\sup_{r \in [h,1]} \tilde{\ell}_*(r) \geq \theta_*\Big\}.
\end{aligned}
$$

Observe that

$$
\begin{aligned}
\tilde{\ell}_*(r) &= \langle \psi_{\tau,\delta_*}, Kf \rangle + \varepsilon \int_{\mathbb{R}^d} \psi_{\tau,\delta_*}(x)\, dW(x) \\
&= \langle f, \varphi_{\tau,\delta_*} \rangle + X(u,r) \\
&= X(u,r) \qquad \text{for } r \in (h,1],
\end{aligned}
$$

because the support of $\varphi_{\tau,\delta_*}$ is $E_\tau$, and $E_\tau \cap G = \varnothing$ for $r \in (h,1]$. Note that, for $\varepsilon$ small enough, $\theta_* \geq c_5 \sigma_X \sqrt{\ln(1/\sigma_X)}$; see (23). Therefore, choosing $\theta = \theta_*$, we have, by Lemma 2,

$$
\begin{aligned}
I_3 &\leq \mathbb{P}\Big\{\sup_{r \in [h,1]} |X(u,r)| \geq \theta_*\Big\} \\
&\leq c_6 \sqrt{\ln\{1/\sigma_X\}}[\sigma_X]^{c_7-1} \leq c_8 \theta_*^{2/\alpha},
\end{aligned}
$$

where the last inequality is obtained by the choice of $C_3$ in (18). Combining bounds for $I_1$, $I_2$ and $I_3$, we complete the proof. □

PROOF OF THEOREM 3. First we note that

$$
\Delta_\infty(\tilde{G}_*, G) \leq \sup_{u \in S^{d-1}} |\tilde{h}_*(u) - h(u)|.
$$

The remainder of the proof is identical to the proof of Theorem 2 with the following minor modifications. When bounding $I_1$ in (28), we take the additional supremum over $u \in S^{d-1}$. Therefore, the expectation $\mathbb{E}\sup_{\tau \in T} |\tilde{\ell}_*(\tau) - \ell_f(\tau)|^2$ should be bounded. This is done exactly as in the proof of Theorem 1, with the only difference being that the stochastic error of the estimator is now considered as the absolute value of the Gaussian process with the index set $T$; then part (ii) of Lemma 1 is used. We note that in this case the constant $C_3$ in (18) should be taken to be larger than in Theorem 1. In the same manner, when bounding $I_3$, the statement (ii) of Lemma 1 is used. □

**Acknowledgment.** The authors would like to thank the referees for their constructive comments and helpful suggestions.

## REFERENCES

[1] ANDERSSEN, R. S. (1980). On the use of linear functionals for Abel-type integral equations in applications. In *The Application and Numerical Solution of Integral Equations* (R. S. Anderssen, F. R. de Hoog and M. A. Lucas, eds.) 195–221. Nijhoff, The Hague. MR0582991
[2] AUBIN, J.-P. (1979). *Applied Functional Analysis*. Wiley, New York. MR0549483
[3] BADDELEY, A. J. (1992). Errors in binary images and an $L^p$ version of the Hausdorff metric. *Nieuw Arch. Wisk. (4)* **10** 157–183. MR1218662
[4] BERTERO, M. and BOCCACCI, P. (1998). *Introduction to Inverse Problems in Imaging*. Institute of Physics, Bristol and Philadelphia. MR1640759
[5] BRANDOLINI, L., RIGOLI, M. and TRAVAGLINI, G. (1998). Average decay of Fourier transforms and geometry of convex sets. *Rev. Mat. Iberoamericana* **14** 519–560. MR1681584
[6] BRUNA, J., NAGEL, A. and WAINGER, S. (1988). Convex hypersurfaces and Fourier transforms. *Ann. of Math. (2)* **127** 333–365. MR0932301
[7] CANDÈS, E. J. and DONOHO, D. L. (2002). Recovering edges in ill-posed inverse problems: Optimality of curvelet frames. *Ann. Statist.* **30** 784–842. MR1922542
[8] DONOHO, D. L. (1999). Wedgelets: Nearly minimax estimation of edges. *Ann. Statist.* **27** 859–897. MR1724034
[9] DONOHO, D. L. and LOW, M. G. (1992). Renormalization exponents and optimal pointwise rates of convergence. *Ann. Statist.* **20** 944–970. MR1165601
[10] EFROMOVICH, S. (1997). Robust and efficient recovery of a signal passed through a filter and then contaminated by non-Gaussian noise. *IEEE Trans. Inform. Theory* **43** 1184–1191. MR1454946
[11] ERMAKOV, M. (1989). Minimax estimation of the solution of an ill-posed convolution type problem. *Problems Inform. Transmission* **25** 191–200. MR1021197
[12] GIHMAN, I. I. and SKOROHOD, A. V. (1974). *The Theory of Stochastic Processes* **1**. Springer, Berlin. MR0346882
[13] GOLBERG, M. (1979). A method of adjoints for solving some ill-posed equations of the first kind. *Appl. Math. Comput.* **5** 123–129. MR0522852
[14] GOLDENSHLUGER, A. (1999). On pointwise adaptive nonparametric deconvolution. *Bernoulli* **5** 907–925. MR1715444
[15] GOLDENSHLUGER, A. and SPOKOINY, V. (2004). On the shape-from-moments problem and recovering edges from noisy Radon data. *Probab. Theory Related Fields* **128** 123–140. MR2027297

[16] GOLDENSHLUGER, A., TSYBAKOV, A. B. and ZEEVI, A. (2006). Optimal change-point estimation from indirect observations. *Ann. Statist.* **34** 350–372.
[17] HALL, P. (1990). Optimal convergence rates in signal recovery. *Ann. Probab.* **18** 887–900. MR1055440
[18] HALL, P. and KOCH, I. (1990). On continuous image models and image analysis in the presence of correlated noise. *Adv. in Appl. Probab.* **22** 332–349. MR1053234
[19] HALL, P., NUSSBAUM, M. and STERN, S. E. (1997). On the estimation of a support curve of indeterminate sharpness. *J. Multivariate Anal.* **62** 204–232. MR1473874
[20] HALL, P. and RAIMONDO, M. (1998). On global performance of approximations to smooth curves using gridded data. *Ann. Statist.* **26** 2206–2217. MR1700228
[21] HÄRDLE, W., PARK, B. U. and TSYBAKOV, A. B. (1995). Estimation of non-sharp support boundaries. *J. Multivariate Anal.* **55** 205–218. MR1370400
[22] KOROSTELËV, A. P. and TSYBAKOV, A. B. (1993). *Minimax Theory of Image Reconstruction. Lecture Notes in Statist.* **82**. Springer, New York. MR1226450
[23] LIGHTHILL, M. J. (1958). *Introduction to Fourier Analysis and Generalised Functions.* Cambridge Univ. Press. MR0092119
[24] LINDENBAUM, M. and BRUCKSTEIN, A. M. (1994). Blind approximation of planar convex sets. *IEEE Trans. Robotics and Automation* **10** 517–529.
[25] MAMMEN, E. and TSYBAKOV, A. (1995). Asymptotical minimax recovery of sets with smooth boundaries. *Ann. Statist.* **23** 502–524. MR1332579
[26] MÜLLER, H.-G. and SONG, K. S. (1994). Maximin estimation of multidimensional boundaries. *J. Multivariate Anal.* **50** 265–281. MR1293046
[27] NEUMANN, M. H. (1997). Optimal change-point estimation in inverse problems. *Scand. J. Statist.* **24** 503–521. MR1615339
[28] RICHARDS, I. and YOUN, H. (1990). *Theory of Distributions*: *A Nontechnical Introduction.* Cambridge Univ. Press. MR1058014
[29] SCHNEIDER, R. (1993). *Convex Bodies*: *The Brunn–Minkowski Theory.* Cambridge Univ. Press. MR1216521
[30] SKIENA, S. S. (1992). Interactive reconstruction via geometric probing. *Proc. IEEE* **80** 1364–1383.
[31] TSYBAKOV, A. (1994). Multidimensional change-point problems and boundary estimation. In *Change-Point Problems* (E. Carlstein, H.-G. Müller and D. Siegmund, eds.) 317–329. IMS, Hayward, CA. MR1477933
[32] VAN DER VAART, A. and WELLNER, J. (1996). *Weak Convergence and Empirical Processes*: *With Applications to Statistics.* Springer, New York. MR1385671

DEPARTMENT OF STATISTICS
HAIFA UNIVERSITY
HAIFA 31905
ISRAEL
E-MAIL: goldensh@stat.haifa.ac.il

GRADUATE SCHOOL OF BUSINESS
COLUMBIA UNIVERSITY
3022 BROADWAY
NEW YORK, NEW YORK 10027
USA
E-MAIL: assaf@gsb.columbia.edu